\documentclass[11pt,a4paper]{article}

\usepackage{amsfonts}
\usepackage{floatflt}
\usepackage{amsthm}

\usepackage{array}

\usepackage{a4wide}

\def\it{\textit} 
\def\bf{\textbf} 
 
\def\ul{\underline}
\def\sf{\textsf} 
 
\def\mc{\mathcal}
\def\mb{\mathbf}
\def\bb{\mathbb}

\def\ds{\displaystyle}

\everymath{\displaystyle}

\newcommand{\be}{\begin{equation}}
\newcommand{\ee}{\end{equation}}
\newcommand{\benum}{\begin{enumerate}}
\newcommand{\eenum}{\end{enumerate}}
\newcommand{\bit}{\begin{itemize}}
\newcommand{\eit}{\end{itemize}}

\newtheorem{thom}{Theorem}
\newtheorem{lemma}[thom]{Lemma}
\newtheorem{rem}[thom]{Remark}

\newtheorem{conj}[thom]{Conjecture}
\newtheorem{defn}[thom]{Definition}

\newtheorem{corol}[thom]{Corollary}

\begin{document}

\title{Computational Complexity, NP Completeness and Optimization
Duality: A Survey}

\author{Prabhu Manyem \\
Department of Mathematics \\
Shanghai University \\
Shanghai 200444, China. \\
Email: \texttt{prabhu.manyem@gmail.com}}


\maketitle

\begin{abstract}
We survey research that studies the connection between the
computational complexity of optimization problems on the one hand, and
the duality gap between the primal and dual optimization problems on the
other.
To our knowledge, this is the first survey that connects the two very
important areas.
We further look at a similar phenomenon in finite model theory
relating to complexity and optimization.
\end{abstract}

\thispagestyle{empty}

\bf{Keywords}.
Computational complexity, Optimization, Mathematical programming,
Duality, Duality Gap, Finite model theory.

\bf{AMS Classification}.
90C25, 90C26, 90C46, 49N15, 65K05, 68Q15, 68Q17, 68Q19 and 68Q25. 


\section{Introduction}\label{sec:intro}

In optimization problems, the \textbf{duality gap} is the difference
between the optimal solution values of the primal problem and the dual
problem.
The relationship between the duality gap and the computational complexity
of optimization problems has been implicitly studied for the last few
decades. 
The connection between the two phenomena has been subtly acknowledged.
The \it{gap} has been exploited to design good approximation algorithms
for NP-hard optimization problems \cite{viggo, hochbaum, vazirani}. 
However, we have been unable to locate a single piece of literature that
addresses this issue explicitly.

This report is an attempt to bring a great deal of evidence
together and specifically address this issue.
Does the existence of polynomial time algorithms for the primal and the
dual problems mean that the duality gap is zero?
Conversely, does the existence of a duality gap imply that either the
primal problem or the dual problem is (or both are) NP-hard?
Is there an inherent connection between computational complexity and
\it{strong duality} (that is, zero duality gap)?

Vecten and Fasbender (independently) were the first to discover the
optimization duality \cite{giannessi2007paper}.
They observed this phenomenon in the \it{Fermat-Torricelli problem}: 
Given a triangle $T_1$, find the equilateral triangle circumscribed
outside $T_1$ with the maximum height.  They showed that this
\it{maximum} height $H$ is equal to the \it{minimum} sum of the distances
from the vertices to $T_1$ to the \it{Torricelli point}\footnote{The
Torricelli point $X$ is indeed the one with the least sum of the distances 
$\vert AX \vert + \vert BX \vert + \vert CX \vert$ from the vertices $A$,
$B$ and $C$ of $T_1$.}.  Thus, this problem enjoys strong duality.

The apparent connection between the duality gap and computational
complexity was considered more than thirty years ago.
Linear Programming (LP) is a well known optimization problem.
In the mid 1970s, before Khachiyan published his ellipsoid algorithm, LP
was thought to be polynomially solvable, precisely because it obeys
\textit{strong duality}; that is, the duality gap is zero.
For a good description of the ellipsoid algorithm, the reader is referred
to the good book by Fang and Puthenpura \cite{fangPuthenpura}.
Strong duality also places the decision version of Linear Programming in
the class NP $\cap$ CoNP; see Lemma \ref{lem:npCoNp} below.

We should stress that this is \bf{not} a survey on Lagrangian duality or
any other form of optimization duality.
Rather, this is a survey on the \bf{connections} and \bf{relationships}
between the computational complexity of optimization problems and
duality.

\section{Definitions}

A few definitions are provided in this section.

\begin{defn}
(\textsf{$\mathsf{D_1 (r)}$: Decision problem corresponding to a given
minimization problem})

{\em Given}.  An objective function $f(\mb{x})$, 
as well as $m_1$ number of constraints 
$\mb{g(x)} = \mb{b}$ and $m_2$ number of constraints
$\mb{h(x)} \ge \mb{c}$,
where $\mb{x} \in \bb{R}^n$ is a vector of variables, 
$\mb{b} \in \bb{R}^{m_1}$ and
$\mb{c} \in \bb{R}^{m_2}$ are constants.  Also given is a parameter
$r \in \bb{R}$.

{\em To Do}.  Determine if the set $\mc{F} = \emptyset$, where
$\mc{F} = \{\mb{x}: ~ \mb{g(x)} = \mb{b}, ~ \mb{h(x)} \ge \mb{c} ~
\mbox{and} ~ f(\mb{x}) \le r\}$.

(Analogously, if the given problem is maximization, then 
\newline
$\mc{F} = \{\mb{x}: ~ \mb{g(x)} = \mb{b}, ~ \mb{h(x)} \ge \mb{c} ~
\mbox{and} ~ f(\mb{x}) \ge r\}$.)
\label{decisionProb}
$\hfill \rule{2.0mm}{2.0mm}$
\end{defn}
$\mc{F}$ is the set of feasible solutions to the decision problem.

\begin{rem}
A word of caution: For decision problems, the term ``feasibility"
includes the constraint on the objective function; if the objective
function constraint is violated, the problem becomes infeasible.
\end{rem}

\begin{defn} 
\cite{bazaraaEtAl}
(\sf{Lagrangian Dual})
\newline
Suppose we are given a minimization problem $P_1$ such as
\begin{equation}
\begin{array}{rl}
 \mbox{Minimize} & f(\mb{x}): X \rightarrow \bb{R} ~~ (X \subseteq
                                        \bb{R}^n), \\ [1.5mm]
 \mbox{subject to} & \mb{g(x)} = \mb{b}, 
        ~~\mb{h(x)} \ge \mb{c}, \\ [1.5mm]
 \mbox{where} & \mb{x} \in \bb{R}^n, 
\mb{b} \in \bb{R}^{m_1} \mbox{ and }
\mb{c} \in \bb{R}^{m_2}.
\end{array}
\label{def:problemP1}
\end{equation}

Let $\mb{u} \in \bb{R}^{m_1}$ and $\mb{v} \in \bb{R}^{m_2}$ be two
vectors of variables with $\mb{v} \ge \mb{0}$.
Let $\ds \mb{e} = (\mb{g}, \mb{h})$.
Assume that 
\bf{b} = $[b_1$ $b_2$ $\cdots$ $b_{m_1}]^T$ and 
\bf{c} = $[c_1$ $c_2$ $\cdots$ $c_{m_2}]^T$ are column vectors.
The feasible region for the primal problem is $X$.

For a given \bf{primal} problem as in $P_1$, the Lagrangian \bf{dual}
problem $P_2$ is defined as follows:
\begin{equation}
\begin{array}{rl}
 \mbox{Maximize} & \theta(\mb{u}, \mb{v}) \\ [1.5mm]
 \mbox{subject to} & \mb{v} \ge \mb{0}, \mbox{where} \\ [1.5mm]
\ds \theta(\mb{u}, \mb{v}) = &
\ds \inf_{\mb{x} \in  \bb{R}^n}
 ~ \{f(\mb{x}) + \sum_{i = 1}^{m_1} u_i (g_i (\mb{x}) - b_i) +
                 \sum_{j = 1}^{m_2} v_j (h_j (\mb{x}) - c_i) \}.
\end{array}
\label{def:problemP2}
\end{equation}
\label{def:LagrangianDual}
$\hfill \rule{2.0mm}{2.0mm}$
\end{defn} 
Note that 
$g_i (\mb{x}) - b_i = 0$ [$h_j (\mb{x}) - c_j \ge 0$] is the $i^{th}$
equality [$j^{th}$ inequality] constraint.

\section{Background: Duality and the classes NP and CoNP}
\label{sec:npCOnp}

(\ul{Note to Reviewers}:  Hope the definitions and explanation in this section
is sufficient.  If not, please let us know where we should expand and
elaborate.)

We now turn our attention to the relationship between the duality of an
optimization problem, and membership in the complexity classes NP and
CoNP of the corresponding set of decision problems.
Decision problems are those with yes/no answers, as opposed to optimization
problems that return an optimal solution (if a feasible solution exists).

Corresponding to $P_1$ defined above in (\ref{def:problemP1}), there is a
set $D_1$ of decision problems,
defined as $D_1 = \{D_1 (r) ~\vert~ r \in \bb{R} \}$. ~The definition of
$D_1 (r)$ was provided in Def. \ref{decisionProb}.

Let us now define the computational classes NP, CoNP and P.
~For more details, the interested reader is referred to either
\cite{viggo} or \cite{papa}.
We begin with the following well known definition:
\begin{defn}
\bf{NP} (respectively \bf{P}) is the class of
decision problems for which there exist non-deterministic (respectively
deterministic) Turing machines which provide Yes/No answers in time that is
polynomial in the size of the input instance.
In particular, for problems in P and NP, if the answer is \bf{yes}, the
Turing machine (TM) is able to provide an ``evidence" (in technical
terms called a {\em certificate}), such as a feasible solution to a
given instance.

The class \bf{CoNP} of decision problems is similar to NP, except for
one key difference: the TM is able to provide a certificate only for
\bf{no} answers.

From the above, it follows that for an instance of a problem in NP $\cap$
CoNP, the corresponding Turing machine can provide a certificate for both
yes and no instances.
\label{npCOnp}
\end{defn}

For example, if $D_1(r)$ in Def. \ref{decisionProb} above is in
NP, the certificate will be a feasible solution; that is, an 
$\mb{x} \in \bb{R}^n$ which obeys the constraints
\begin{equation}
\mb{g(x)} = \mb{b}, ~ 
\mb{h(x)} \le \mb{c} \mbox{ and }
f(\mb{x}) \ge r.
\label{allConstraints}
\end{equation}
On the other hand, if $D_1(r) \in$ CoNP, the certificate will be an
$\mb{x} \in \bb{R}^n$ that violates at least one of the $m_1 + m_2 + 1$
constraints in (\ref{allConstraints}).

\begin{rem}
For problems in NP, for Yes instances, extracting a solution from the
certificate is not always an efficient (polynomial time) task.
Similarly, in the case of CoNP, pinpointing a violation from a Turing
machine certificate
\footnote{We thank WenXun Xing and PingKe Li of Tsinghua University
(Beijing) for ponting out the above.}
is not guaranteed to be efficient either.
\end{rem}

\begin{rem}
P $\subseteq$ NP, because any computation that can be carried out by a
deterministic TM can also be carried out by a non-deterministic TM.
The problems in P are decidable deterministically in polynomial time.

The class P is the same as its complement Co-P.
~That is, P is closed under complementation.

Furthermore,
Co-P ($\equiv$ P) is a subset of CoNP.
~We know that P is a subset of NP.
~Hence $P \subseteq$ NP $\cap$ CoNP.
Thus for an instance of a problem in P, the corresponding Turing
machine can provide a certificate for both yes and no instances.
\label{PsubsetConp}
\end{rem}

We are now ready to define what is meant by a \it{tight dual}, and how it
relates to the intersection class of problems, NP $\cap$ CoNP.
~Note that for two problems to be tight duals, it is sufficient if they are
tight with respect to just one type of duality (such as Lagrangian duality,
for example).

\begin{defn}
\textsf{Tight duals and the class TD}.
Two optimization problems $P_a$ and $P_b$ are \bf{dual to
each other} if the dual of one problem is the other.

Suppose $P_a$ and $P_b$ are dual to each other, with zero duality gap;
then we say that $P_a$ and $P_b$ are \bf{tight duals}.
For any $r \in \bb{R}$, let $D_a (r)$ and $D_b (r)$ be the decision
versions of $P_a$ and $P_b$ respectively.

Let \bf{TD} be the class of all decision problems whose optimization
versions have tight duals.  That is, \bf{TD} is the set of all problems 
$D_a (r)$ and $D_b (r)$ for any $r \in \bb{R}$.
\label{def:tightDuals}
\end{defn}

\begin{rem}
{A word of caution}.  Tight duality is not the same as strong duality.
For strong duality, it is sufficient if there exist feasible solutions
to the primal and the dual such that the duality gap is zero.
For tight duality to hold, we also require that the primal problem $P_a$
and the dual problem $P_b$ be dual to one another.
\end{rem}

One way in which duality gaps are related to the classes NP and CoNP is
as follows:
\begin{lemma}
\cite{papa}
TD $\subseteq$ NP $\cap$ CoNP.
\label{lem:npCoNp}
\end{lemma}

From Remark \ref{PsubsetConp} and Lemma \ref{lem:npCoNp}, we know that
both TD and P are subsets of NP $\cap$ CoNP.~
But is there a containment relationship between TD and P?~
That is, is either TD $\subseteq$ P or P $\subseteq$ TD?~
This is the subject of further study in this paper, with particular
reference to Lagrangian duality.
\begin{rem}
We should mention that in several cases, given a primal problem $P$, even
if we are able to find a dual problem $D$ such that the dual of $P$ is
$D$, it \it{does not} necessarily follow that the dual of $D$ is $P$.
~That is,the dual of the dual \it{need not} be the primal.
$P$ and $D$ are \it{not necessarily duals of each other}.
We do not include such $(P, D)$ pairs in \bf{TD}.
~Among primal-dual pairs of problems, TD is a restricted class.
\label{rem:TDpairs}
\end{rem}

\section{Lack of Strong Duality results in NP hardness}
\label{sec:gapMeansHard}

In this section, we will review results from the literature, which show
that the lack of strong duality imply that the optimization problem in
question is NP-hard, assuming that the primal problem obeys the
\it{constraint qualification} assumption as stated below in Def.
\ref{convexProg}.
Here we work with Lagrangian duality.
Results for other types of duality such as Fenchel, geometric and canonical
dualities require further investigation.

Let us define what we mean by weak duality (as opposed to tight duality
and strong duality):

\begin{defn}
Given a primal problem $P_1$ and a dual problem $P_2$, as defined in
Def. \ref{def:LagrangianDual}, the pair ($P_1, P_2$) is said to
obey \bf{weak duality} if $\theta(\mb{u}, \mb{v}) \le f(\mb{x})$, for
every feasible solution \bf{x} to the primal and every feasible solution
(\bf{u, v}) to the dual.
\label{def:weakDuals}
$\hfill \rule{2.0mm}{2.0mm}$
\end{defn}
\begin{defn}
In Def. (\ref{def:weakDuals}), if
\newline
(i) the inequality is replaced by an equality, and 
\newline
(ii) there exist a primal feasible solution $\mb{\bar{x}}$ and a dual
feasible solution $(\mb{\bar{u}}, \mb{\bar{v}})$ such that the
equality in (i) holds,
\newline
then the pair ($P_1, P_2$) is said to be obey \bf{strong duality}.
\end{defn}

The following theorem from \cite{bazaraaEtAl} guarantees that the 
feasible solutions to Lagrangian dual problems (\ref{def:problemP1}) and
(\ref{def:problemP2}) indeed obey \it{weak duality}:
\begin{thom}
If \bf{x} is a feasible solution to the primal problem in
(\ref{def:problemP1}) and (\bf{u}, \bf{v}) is a feasible solution to the
dual problem in (\ref{def:problemP2}), then 
$f(\mb{x}) \ge \theta(\mb{u}, \mb{v})$.
\label{weakDuality}
$\hfill \rule{2.0mm}{2.0mm}$
\end{thom}

We shall now define a special type of convex program, called a 
\it{convex program with constraint qualification}, which is one with an
assumption about the existence of a feasible solution in the interior of
the domain.
\begin{defn}
\textsf{Convex program (convex optimization problem)}.
\newline
{\em Given}. A convex set $X \subset \bb{R}^n$, two convex functions
$f(\mb{x}): \bb{R}^n \rightarrow \bb{R}$
and $\mb{g(x)}: \bb{R}^n \rightarrow \bb{R}^{m_1}$, 
\newline
as well as an affine function $\mb{h(x)}:
\bb{R}^n \rightarrow \bb{R}^{m_2}$.

{\em To do}.
Minimize $f(\mb{x})$, subject to 
$\mb{g(x) \le 0}$, 
~$\mb{h(x) = 0}$ and 
$\mb{x} \in X$.
\label{convexProg}
\end{defn}

\begin{defn}
\textsf{Convex program with constraint qualification}.
\newline
Same as the optimization problem in Def. \ref{convexProg}, except
that we include the following 
{\em constraint qualification assumption}.
There is an $\mb{x_0} \in X$ such that
$\mb{g(x_0) < 0}$,
~$\mb{h(x_0) = 0}$, and
~$\ds \mb{0} \in \mbox{int} ~ \mb{h}(X)$,
where 
$\ds \mb{h}(X) = \bigcap_{\mb{x} \in X} h(\mb{x})$.
\label{convexProgWithCC}
\end{defn}

(\it{Note}:
Of course, the functions above can be written in the same form as in Def.
\ref{def:problemP1} and \ref{def:problemP2}.  In such a case, we can
define $(\mb{g(x) - b})$ to be a convex function and ~$(\mb{h(x) -
c})$ to be an affine function, where $\mb{b} \in \bb{R}^{m_1}$ and 
$\mb{c} \in \bb{R}^{m_2}$.)

\begin{defn}
In Def. \ref{convexProg}, if any of the ($m+1$) functions 
$f(\mb{x}): \bb{R}^n \rightarrow \bb{R}$
and $\mb{g(x)}: \bb{R}^n \rightarrow \bb{R}^{m_1}$ is not convex, then
the optimization problem is said to \bf{non-convex}.
\label{def:nonConvexProg}
\end{defn}

For the remainder of this section, we will assume primal constraint
qualification; that is, we assume that constraint qualification is applied
to the primal optimization problem.  The following theorem provides
sufficient conditions under which strong duality can occur:

\begin{thom}
\bf{Strong Duality} \cite{bazaraaEtAl}.
If (i) the primal problem is given as in Def. \ref{convexProg}, and
(ii) the primal and dual problems have feasible solutions, 
then the primal and dual optimal solution
values are equal (that is, the duality gap is zero):
\begin{equation}
\begin{array}{rl}
& \inf
  \{f(\mb{x}): \mb{x} \in X, ~ \mb{g(x) \le 0},
~ \mb{h(x) = 0} \}
=
\sup \{\theta(\mb{u}, \mb{v}): \mb{v \ge 0} \}, \\ [1.5mm]
& \theta(\mb{u}, \mb{v}) = 
\inf_{\mb{x} \in  X}
  \{f(\mb{x}) + \sum_{i = 1}^{m_1} u_i g_i (\mb{x}) +
                 \sum_{j = 1}^{m_2} v_j h_j (\mb{x}) \},
\end{array}
\end{equation}
where 
$\ds \theta(\mb{u}, \mb{v})$ is the dual objective function.
\label{thom:strongDuality}
\end{thom}

Using the contrapositive statement of Theorem \ref{thom:strongDuality}, we
get the following result:

\begin{corol}
(to Theorem \ref{thom:strongDuality})
If there exists a duality gap using
Lagrangian duals, then either the primal or the dual is not a convex
optimization problem. (Remember, we are assuming constraint qualification.)
\label{corolToStrongDuality}
\end{corol}

The \it{Subset Sum} problem is defined as follows: Given a set $S$ of positive
integers $\{d_1, d_2, \cdots, d_k\}$ and another positive integer $d_0$, is
there a subset $P$ of $S$, such that the sum of the integers in $P$ equals
$d_0$?

Using a polynomial time reduction from the \it{Subset Sum} problem to a
non-convex optimization problem (see Def. \ref{def:nonConvexProg}),
Murty and Kabadi (1987) showed the following:
\begin{thom}
\cite{murtyKabadi}
If an optimization problem is non-convex, it is NP-hard.
\label{thm:kabadiMurty}
\end{thom}
In certain cases, non-convex problems have an equivalent convex
formulation, for example, through strong duality.
Such a dual transformation, where a convex problem $B$ is a dual
of a non-convex problem $A$ such that the duality gap between them is
zero, is called \textsf{hidden convexity} \cite{beck2009,
benTalTeboulle}.
In such cases, the reformulated convex problem is also NP-hard; otherwise
the primal non-convex problem can be solved efficiently, thus violating
Theorem \ref{thm:kabadiMurty}.

\begin{defn}
\textsf{Standard Quadratic Program (SQP)} \cite{bomzeDeClerk2002}.
Minimize the function~ $\ds \mb{x}^T Q \mb{x}$,
where $\mb{x} \in \Delta$  (the standard simplex in
$\mathbb{R}^n_+$). ~The $n$ vertices of $\Delta$ are at a unit distance 
(in the positive direction) along each of the $n$ axes of $\bb{R}^n$.
~$Q$ is a given symmetric matrix in $\bb{R}^{n \times n}$.
\end{defn}

The converse of Theorem \ref{thm:kabadiMurty} is not true. A convex
optimization problem in general is NP-hard, SQP being an example.
~SQP is non-convex; however in \cite{bomzeDeClerk2002}, Bomze and de
Clerk prove that it has an exact convex reformulation as a \it{copositive
programming} problem. 
SQP is known to be NP-hard, since its decision version contains the
max-clique problem in graphs as a special case.
From this, it follows that {copositive programming} is also NP-hard; see
\cite{mirjamDur} for more on this topic.

From Corollary
\ref{corolToStrongDuality} and
Theorem \ref{thm:kabadiMurty}, it follows that
\begin{thom}
Assuming constraint qualification, 
if there exists a duality gap using Lagrangian duals, then either the
primal or the dual is NP-hard.
\end{thom}

These results are true for Lagrangian duality.
For other types of duality such as Fenchel, geometric and canonical
dualities, this requires further investigation.

\section{Does Strong Duality Imply Polynomial Time Solvability?}
\label{sec:noGapMeansEasy}

At this time, such a proof (of whether a duality gap of zero implies
polynomial time solvability of the primal and the dual problems) appears
possible only for very simple problems, since estimating the duality gap
appears extremely challenging for many problems.

Some of the problems where this is true include Convex Programming (and in
particular, Linear Programming), where both the primal and the dual
problems are convex; the polynomial time algorithms are derived from
interior point methods.
For example, see the book by Nesterov and Nemirovskii \cite{nn94}.

In a recent paper \cite{manyem10}, we have demonstrated an additional
computational benefit arising from strong duality. Primal-dual
problem pairs that are polynomially solvable and obey strong duality can
be solved by a single call to a \it{decision Turing machine}, that is, a
Turing machine that provides a Yes/No answer (if the answer is yes, then
it can provide a feasible solution which supports the Yes answer).
Previously, it was only known that optimization problems require multiple
calls to a decision Turing machine (for example, doing a binary search on
the solution value to obtain an optimal solution).
For more details, the reader is referred to \cite{manyem10}.

More investigation is needed to answer the question
\it{Does Strong Duality Imply Polynomial Time Solvability?
} in a general setting.

However, there have been some results recently using Canonical duality for
certain types of quadratic programs \cite{gao2004, xingFangGaoSheuZhang}.
Consider a standard quadratic programming (primal) problem:
\begin{equation}
\begin{array}{rl}
\mbox{Minimize} & P(\mb{x}) = \frac{1}{2}\mb{x}^T A \mb{x}
                    + \mb{b}_0^T \mb{x} \\ [1.5mm]
 \mbox{subject to} & \mb{b}_i^T \mb{x} \le {c_i}, ~ 1 \le i \le m,
\end{array}
\label{canonicalProb}
\end{equation}
where $A = A^T \in \bb{R}^{n \times n}$, 
$\mb{b}_i \in \bb{R}^n$ for $1 \le i \le m$,
and $\mb{x} \in \bb{R}^n$.
~$A$ is a symmetric matrix.
We can write the Lagrangian function as
\begin{equation}
L(\mb{x}, \mb{\lambda}) = \frac{1}{2}\mb{x}^T A \mb{x} + 
    \left( \mb{b}_0 + \sum_{i=1}^m \lambda_i \mb{b}_i \right)^T \mb{x} - 
      \sum_{i=1}^m \lambda_i c_i, ~ \lambda_i \ge 0.
\end{equation}
The first order necessary condition among the Karush-Kuhn-Tucker (KKT)
conditions yields
$\ds A\mb{x} + \mb{b}_0 + \sum_{i=1}^m \lambda_i \mb{b}_i = 0$, from which
we get the value of $\mb{x}$ as
\begin{equation}
\mb{x} = -A^{-1} \left( \mb{b}_0 + \sum_{i=1}^m \lambda_i \mb{b}_i \right),
\end{equation}
assuming that $A$ is an invertible matrix.  Substituting this value of
\bf{x} back into the Lagrangian function yields
\begin{equation}
Q(\lambda) = -\frac{1}{2} \left( \mb{b}_0 + \sum_{i=1}^m \lambda_i \mb{b}_i
\right)^T A^{-1} \left( \mb{b}_0 + \sum_{i=1}^m \lambda_i \mb{b}_i \right)
- \sum_{i=1}^m \lambda_i c_i.
\end{equation}
We then get a dual problem, the canonical dual $P^d$, as follows:
\begin{equation}
\begin{array}{rl}
 \mbox{Maximize} & Q(\lambda), \\ [1.5mm]
 \mbox{subject to} & A \succ 0 ~ \mbox{and} ~ \lambda \ge 0.
\end{array}
\end{equation}

(The relation $A \succ 0$ means the matrix $A$ should be positive definite.)

The positive definiteness of $A$ has been found to be a sufficient
condition.  Whether it is a necessary condition is not known yet.

Now, $Q(\lambda)$ is a concave function, to be maximized in the dual
variable $\lambda$;
hence $P^d$ can be solved efficiently, since we are
maximizing a concave function.
Thus in general, we can get a lower bound for the primal problem quickly;
and in some cases, we can get a strong dual with zero duality gap, which
provides an optimal solution for the primal problem in polynomial time.

Canonical duality was first developed to address the problem of (possibly
large) duality gaps in Lagrangian duality in the context of problems in
analytical mechanics \cite{gaoRuanSherali2009}.
Furthermore, large duality gaps can also be found while solving
{non-convex} optimization problems using {Fenchel-Morrow-Rockafellar}
duality.  The so-called \it{canonical dual transformation} can be used to
formulate a strong-dual problem (that is, one with a zero duality gap).
Quoting from \cite{gaoRuanSherali2009}, \it{the primal problem can be
made equal to its canonical dual in the sense that they have the same KKT
points}.

\subsection{Preliminary Results:
Canonical Duality and the Complexity Classes NP and CoNP}
\label{sec:optProblems}

We consider quadratic programming
\footnote{Thanks to Shu-Cherng Fang (NCSU, USA) for his input here.} 
problems with a single quadratic
constraint (QCQP) \cite{xingFangGaoSheuZhang}.
The primal problem $P_0$ is given by:
\begin{equation}
\begin{array}{rl}
 \mbox{Minimize} & P(\mb{x}) = \frac{1}{2}\mb{x}^T A \mb{x}
                    - \mb{f}^T \mb{x} \\ [1.5mm]
 \mbox{subject to} & \frac{1}{2} \mb{x}^T B \mb{x} \le {\mu}, 
\end{array}
\label{eq:quadConstraint}
\end{equation}
where $A$ and $B$ are non-zero $n \times n$ symmetric matrices, 
$\mb{f} \in \mb{R}^n$, and $\mu \in \mb{R}$.
($A$, $B$, $\mb{f}$ and $\mu$ are given.)
Corresponding to the primal $P_0$, the canonical dual problem $P^d$ is as
follows:
\begin{equation}
\begin{array}{rl}
 \sup_{\sigma} & P^d (\sigma) = - \frac{1}{2} 
                   \mb{f}^T (A + \sigma B)^{-1} \mb{f}
                    - \mu \sigma \\ [2mm]
 \mbox{subject to} & \sigma \in \mc{F} = \{\sigma \ge 0 ~\vert ~
                     A + \sigma B \succ 0 \},
\end{array}
\label{eq:quadConstraintDual}
\end{equation}
assuming that $\mc{F}$ is non-empty.
To ensure that $(A + \sigma B)$ is invertible, we have assumed that $A$
and $B$ are symmetric and non-zero.
Again, the positive definiteness of $(A + \sigma B)$ for the dual problem
has been found to be a sufficient condition, and it is not known yet
whether it is a necessary condition.

The following theorem appeared in 
\cite{xingFangGaoSheuZhang}:
\begin{thom}
(Strong duality theorem)
When the maximum value $P^d (\sigma^{\ast})$ in
(\ref{eq:quadConstraintDual}) is finite, strong duality between the primal
problem $P_0$ and the dual problem $P^d$ holds, and the optimal solution
for $P_0$ is given by $\ds x^{\ast} = (A + \sigma^{\ast} B)^{-1} \mb{f}$.
\label{thm:qcqpStrongDual}
\end{thom}
However, the pair ($P_0$, $P^d$) may not be tight duals (i.e., they may
not belong to the class $TD$).
~To show that ($P_0$, $P^d$) $\in TD$, it also needs to be shown that the
dual of $P^d$ is $P_0$. (Recall Remark \ref{rem:TDpairs}.)

\begin{rem}
We have not provided a detailed analysis of how the canonical problem can
be derived.  The interested reader is referred to 
\cite{xingFangGaoSheuZhang}.  Our goal is to illustrate an instance of a
pair of primal and dual problems that obey strong duality, and where an
optimal solution can be obtained efficiently (in polynomial time).
\end{rem}

Strong duality for a variation of QCQP was established by Mor{\'e} and
Sorensen in 1983 \cite{moreSorensen83} (as described in
\cite{Wolkowicz2000}); they call theirs as the Trust Region Problem (TRP).
The problem is to minimize a quadratic function $q_0$, subject to a norm
constraint $\ds \mb{x}^T\mb{x} \le \delta^2$.
This is an example of a non-convex problem where strong duality holds.
Rendl and Wolkowicz \cite{rendlWolko1997} show that the TRP can be
reformulated as an unconstrained concave maximization problem, and hence
can be solved in polynomial time by using the interior point method
\cite{nn94}.

Recall that decision problems are those with Yes/No answers.
Combining Theorem \ref{thm:qcqpStrongDual}
with Lemma \ref{lem:npCoNp} above, we can conclude
that 
\begin{rem}
Decision versions of Problems $P_0$ in (\ref{eq:quadConstraint}) and
$P^d$ in (\ref{eq:quadConstraintDual}) are members of the intersection
class NP $\cap$ CoNP.
\end{rem}

Furthermore, from Theorem 3 in \cite{xingFangGaoSheuZhang}, it is easy to
see that problems $P_0$ and $P^d$ can be solved in polynomial time.  The
authors in  \cite{xingFangGaoSheuZhang} use what they call a
``boundarification" technique which moves the analytical solution
$\bar{x}$ to a global minimizer $x^{\ast}$ on the boundary of the primal
domain.
This strengthens the conjecture that 
\begin{conj}
If two optimization problems $P_a$ and $P_b$ are such that one of them is
the dual of the other, that is, they exhibit strong duality, 
then both are polynomially solvable.
\end{conj}

As hinted earlier, the canonical dual feasible space $\cal F$ in
(\ref{eq:quadConstraintDual}) may be empty.
In such a case, strong duality still holds on a feasible domain in which
the matrix $(A + \sigma B)$ is not definite.
However, the question as to how to solve the canonical dual problem is
still open.  It is conjectured in \cite{gao2007} that the primal problem
(from which the canonical dual problem is derived) could be NP-hard.

The observations above are for quadratic programming problems with a
single quadratic constraint.
It would be interesting to see what happens for quadratic programming
problems with two constraints, whether strong duality still holds, and
whether both the primal and the dual are still polynomially solvable.

\it{Semidefinite Programming (SDP)}.
Ramana \cite{ramana} exhibited strong duality for the SDP problem.
However, the complexity of SDP is unknown; it was shown in \cite{ramana}
that the decision version of SDP is NP-complete if and only if NP = CoNP.

(\ul{Note}:  There have been some published papers which say
that SDP is polynomial time solvable.  This is NOT correct, as the above
result in \cite{ramana} shows.)

\section{Descriptive Complexity and Fixed Points}

On a final note, we would like to briefly describe a similar phenomenon
which occurs in the field of Descriptive Complexity, which is the
application of Finite Model theory to computational complexity.
In particular, we would like to mention least fixed point (LFP)
computation.
A full description would be beyond the scope of this paper.  However, we
would like to briefly mention a few related concepts and phenomena.

For a good description of least fixed points (LFP) in existential second
order (ESO) logic, the reader is referred to \cite{gradelAnd7others}
(chapters 2 and 3) and \cite{EF99}.
If the input structures are ordered, then expressions in LFP logic can
describe polynomial time (PTIME) computation \cite{gradelAnd7others}.

The input instance to an LFP computation consists of a \it{structure}
\bf{A}, which includes a domain set $A$ and a set of (first order)
relations $R_i$, each with arity $r_i$, $1 \le i \le J$.
The LFP computation works by a stagewise addition of tuples from $A$, to
a new relation $P$ (of some arity $k$).
If $P_i$ represents the relation (set of tuples) after stage $i$,
then $P_i \subseteq P_{i+1}$.
The transition from $P_i$ to $P_{i+1}$ is through an operator $\Phi$,
such that $P_{i+1} = \Phi(P_i)$.
At the beginning, $P$ is empty, that is, $P_0 = \emptyset$.
For some value of $i$, say when $i = f$, if $P_f = P_{f+1}$, a \it{fixed
point} has been reached.

Without going into details, let us just say that such a fixed point,
reached as above, is also a \it{least} fixed point (LFP) if the operator
$\Phi$ can be chosen in a particular manner.
~The interested reader is referred to \cite{gradelAnd7others} (chapter 2)
for details.

Note that the number of elements in $P$ can be at most $|A|^k$ (where
$|A|$ is the number of elements in $A$), which is polynomial in the size
of the domain.
Hence $f \le |A|^k$, so an LFP is achieved within a polynomial number of
stages.

Similar to LFP, we can also define a \it{greatest} fixed point (GFP).
~This is obtained by doing the reverse; we start with the entire set
$A^k$ of $k$-ary tuples from the universe $A$, and then removing tuples
from $P$ in stages.
At the beginning, $P_0 = A^k$.
~In further stages, $P_i \supset P_{i+1}$.
~The GFP is reached at stage $g$ if $P_g = P_{g+1}$.

The logic that includes LFP and GFP expressions is known as \bf{LFP
logic}. It expresses decision problems (those with a Yes/No answer), such
as those in Def. \ref{decisionProb} and \ref{npCOnp}.
To be feasible, a solution should also obey the objective function
constraint ($f(\mb{x}) \ge K$ or $f(\mb{x}) \le K$).

The LFP computation expresses decision problems based on maximization.
Before the fixed point is reached, the solution is infeasible; that is,
the number of tuples in the fixed point relation $P$ is insufficient.
However, once the fixed point is reached, the solution becomes feasible.
Similarly, the GFP computation expresses decision problems based on
minimization.

\bf{Problem}.
An interesting problem arising in LFP Logic is this:
For what type of primal-dual optimization problem pairs will the LFP and
GFP computation meet at the same fixed point?  Does this mean that such a
pair is polynomially solvable?

\section{Conclusion and Further Study}
Let us again stress that this is \bf{not} a survey on Lagrangian
duality or any other form of optimization duality.
Rather, this is a survey on the \bf{connections} and \bf{relationships}
between the computational complexity of optimization problems and duality.

In this paper, we have touched the tip of the iceberg on a very
interesting problem, that of connecting the computational hardness of an
optimization problem with its duality characteristics.
A lot more study is required in this area.

Another issue is that of \it{saddle point} for Lagrangian duals.
This is a decidable problem; we can do brute force and find the
primal and dual optimal solutions; this will tell us if there is a
duality gap.  If the gap is zero, then there is a saddle point.

(Jeroslow \cite{jeroslow} showed that the integer programming problem
with quadratic constraints is undecidable if the number of variables is
unbounded, which is an extreme condition.  However, if each variable
has a finite upper and lower bound, then the number of solutions is
finite and thus it is possible to determine the best solution in finite
time.)

However, this problem would be NP-complete, unless we can tell whether it
has a saddle point by looking at the structure of the problem or by
running a polynomial time algorithm.

We hope that this paper will motivate further research in this very
interesting topic.

\section{Acknowledgments}

This research was supported by 
(i) a visiting fellowship from the National Cheng Kung University (NCKU),
Tainan, Taiwan;
(ii) the Shanghai Leading Academic Discipline Project \#S30104; and
(iii) the National Natural Science Foundation of China (\# 11071158).
~Support from all sources is gratefully acknowledged.


\end{document}